\documentclass[12pt,a4paper,fleqn]{article}%
\usepackage{a4wide,amsfonts,amsmath,latexsym,amssymb,euscript,graphicx,units,mathrsfs}
\usepackage{graphicx}
\usepackage{color}
\usepackage{amssymb}
\usepackage{amssymb}
\usepackage[T1]{fontenc}
\usepackage{latexsym}
\usepackage{xypic}
\usepackage{eufrak}
\usepackage{euscript}
\usepackage{amsmath}
\usepackage{verbatim}
\usepackage{fancyhdr}
\usepackage[english]{babel}
\usepackage{mathrsfs}
\usepackage{units}
\usepackage{amsfonts}%
\newtheorem{prop}{Proposition}[section]
\newtheorem{cor}[prop]{Corollary}
\newtheorem{lemme}[prop]{Lemma}
\newtheorem{rem}[prop]{Remark}
\newtheorem{thm}[prop]{Theorem}

\renewcommand{\geq}{\geqslant}

\newcommand{\R}{\mathbb{R}}

\def\HH{\EuFrak H}

\def\1{{\mathbf{1}}}

\def\1{{\mathbf{1}}}
\def\0.5{{\frac{1}{2}}}

\newcommand{\qed}{\nopagebreak\hspace*{\fill}
{\vrule width6pt height6ptdepth0pt}\par}

\begin{document}

\title{Density estimates and concentration inequalities\newline with Malliavin
calculus }
\author{Ivan Nourdin\thanks{Laboratoire de Probabilit{\'{e}}s et Mod{\`{e}}les
Al{\'{e}}atoires, Universit{\'{e}} Pierre et Marie Curie, Bo{\^{\i}}te
courrier 188, 4 Place Jussieu, 75252 Paris Cedex 5, France,
\texttt{ivan.nourdin@upmc.fr}}\ \ and\ \ Frederi G. Viens\thanks{Dept.
Statistics and Dept. Mathematics, Purdue University, 150 N. University St.,
West Lafayette, IN 47907-2067, USA, \texttt{viens@purdue.edu}}
\and \textit{Universit\'{e} Paris 6}\thinspace\ and \thinspace\textit{Purdue
University}}
\date{}
\maketitle

\begin{abstract}
We show how to use the Malliavin calculus to obtain density estimates of the
law of general centered random variables. In particular, under a
non-degeneracy condition, we prove and use a new formula for the density
$\rho$\ of a random variable $Z$ which is measurable and differentiable with
respect to a given isonormal Gaussian process.
%%J'ai un peu changé (je n'aime pas trop les formules dans l'abstract). Voici l'ancienne version:
%%an abstract Wiener space:
%%\[
%%{\small \rho(z)=}\frac{E|Z|}{2g(z)  }{\small \mathrm{\exp}(}%
%%\int_{0}^{z}\frac{x\,dx}{g\left(  x\right)  }{\small ),}%
%%\]
%%where $g(z)=E(<DZ,-DL^{{\small -1}}Z>|Z=z)$, where $D$\ is the Malliavin
%%derivative, and where $L$\ is the Ornstein-Uhlenbeck operator. 
Among other results, we apply our techniques to bound the density of the maximum of a
general Gaussian process from above and below; several new results ensue,
including improvements on the so-called Borell-Sudakov inequality. We then
explain what can be done when one is only interested in or capable of deriving
concentration inequalities, i.e. tail bounds from above or below but not
necessarily both simultaneously.

\end{abstract}

\bigskip

{\small \noindent\textbf{Key words:} Malliavin calculus; density estimates;
concentration inequalities; fractional Brownian motion;
Borell-Sudakov inequality; suprema of Gaussian processes.~\newline}

{\small \noindent\textbf{2000 Mathematics Subject Classification:} 60G15;
60H07.~\newline}

\section{Introduction}

Let $N$ be a zero-mean Gaussian random vector, with covariance matrix
$K\in\mathscr{S}_{n}^{+}(\mathbb{R})$. Set $\sigma_{\max}^{2}:=\max_{i}K_{ii}%
$, and consider
\begin{equation}
Z=\max_{1\leqslant i\leqslant n}N_{i}-E\big(\max_{1\leqslant i\leqslant
n}N_{i}\big). \label{def-z}%
\end{equation}
It is well-known, see e.g. Vitale \cite{vitale}, that for all $z>0$,%
\begin{equation}
P\big(Z\geqslant z\big)\leqslant\mathrm{exp}\left(  -\frac{z^{2}}%
{2\,\sigma_{\max}^{2}}\right)  \quad\mbox{if $z>0$.} \label{vitale-b}%
\end{equation}
{{The corresponding left-tail probability bound analogue of (\ref{vitale-b})
also holds, see e.g. Borell \cite{borell}:}}
\begin{equation}
P\big(Z\leqslant-z\big)\leqslant\mathrm{exp}\left(  -\frac{z^{2}}%
{2\,\sigma_{\max}^{2}}\right)  \quad\mbox{if $z<0$.} \label{borell-b}%
\end{equation}
Of course, we can combine (\ref{vitale-b}) and (\ref{borell-b}) to get,%
\begin{equation}
P\big(|Z|\geqslant z\big)\leqslant2\,\mathrm{exp}\left(  -\frac{z^{2}%
}{2\,\sigma_{\max}^{2}}\right)  \quad\mbox{if $z>0$.} \label{vitale-borell-b}%
\end{equation}
{{Inequality (\ref{vitale-borell-b}) is a special case of bounds for more
general Gaussian fields. Such bounds are often collectively known as
Borell-Sudakov inequalities. These can be extended much beyond the Gaussian
realm; see for instance the book of Ledoux and Talagrand \cite{LT}. Yet these
Borell-Sudakov inequalities can still be improved, even in the Gaussian
framework; this is one of the things we will illustrate in this paper.}%
}\vspace{0.1in}

{{Inequality (\ref{vitale-borell-b}) is also a special case of results based
on almost sure bounds on a random field's Malliavin derivatives, see Viens and
Vizcarra \cite{VVJFA}. }}While {{that paper }}uncovered a new way to relate
scales of regularity and fractional exponential moment conditions with
iterated Malliavin derivatives, it failed to realize how best to use these
derivatives when seeking basic estimates such as (\ref{vitale-borell-b}). In
the present paper, our aim is to explain how to use Malliavin calculus more
efficiently than in \cite{VVJFA} in order to obtain bounds like
(\ref{vitale-b}) or (\ref{borell-b}), and even often much better. For
instance, by applying our machinery to $Z$ defined by (\ref{def-z}), we obtain
the following.

\begin{prop}
\label{rhointroprop}With $N$ and $Z$ as above, if $\sigma_{\min}^{2}%
:=\min_{i,j}K_{ij}>0$, with $\sigma_{\max}^{2}:={{\max_{i}K_{ii}}}$, the
\textit{density} $\rho$ of $Z$ exists and satisfies, for almost all
$z\in\mathbb{R}$,
\begin{equation}
\frac{E|Z|}{2\sigma_{\max}^{2}}\,\mathrm{exp}\left(  -\frac{z^{2}}%
{2\,\sigma_{\min}^{2}}\right)  \leqslant\rho(z)\leqslant\frac{E|Z|}%
{2\sigma_{\min}^{2}}\,\mathrm{exp}\left(  -\frac{z^{2}}{2\,\sigma_{\max}^{2}%
}\right)  . \label{rho-intro}%
\end{equation}

\end{prop}

This proposition generalizes immediately (see Proposition \ref{max-intro} in
Section \ref{density} below) to the case of processes {{defined on an interval
$[a,b]\subset\mathbb{R}$}}. To our knowledge, that result is the first
instance where the density of the maximum of a general Gaussian process is
estimated from above and below. {{As an explicit application, let us mention
the following result, concerning the centered maximum of a fractional Brownian
motion (fBm), which is proved at the end of Section \ref{density}.}}

\begin{prop}
\label{pr} Let $b>a>0$, and $B=(B_{t},\,t\geqslant0)$ be a fractional Brownian
motion with Hurst index $H\in(1/2,1)$. Then the random variable $Z=\sup
_{[a,b]}B-E\big(\sup_{[a,b]}B\big)$ has a density $\rho$ satisfying, for
almost all $z\in\mathbb{R}$:
\begin{equation}
\frac{E|Z|}{2b^{2H}}\,e^{-\frac{z^{2}}{2\,a^{2H}}}\leqslant\rho(z)\leqslant
\frac{E|Z|}{2a^{2H}}\,e^{-\frac{z^{2}}{2\,b^{2H}}}. \label{ffbbmm}%
\end{equation}

\end{prop}

{Of course, the interest of this result lies in the fact that the exact
distribution of $\sup_{[a,b]}B$ is still an open problem when $H\neq1/2$.
Moreover, note that introducing a degeneracy in the covariances for stochastic
processes such as fBm has dire consequences on their supremas' tails; for
instance, with $a=0$, $Z$ has no left hand tail, since $Z\geqslant
-E\big(\sup_{[0,b]}B\big)$ a.s., and therefore $\rho$ is zero for $z$ small
enough.\vspace{0.1in}}

Density estimates of the type (\ref{rho-intro}) may be used immediately to
derive tail estimates by combining simple integration with the following
classical inequalities:
\[
\frac{z}{1+z^{2}}~e^{-\frac{z^{2}}{2}}\leqslant\int_{z}^{\infty}%
e^{-\frac{y^{2}}{2}}dy\leqslant\frac{1}{z}~e^{-\frac{z^{2}}{2}}\quad
\mbox{for all
$z>0$.}
\]
The two tails of the supremum of a Gaussian vector or process are typically
not symmetric, and neither are the methods for estimating them; this poses a
problem for the techniques used in \cite{vitale} and \cite{borell}, and for
ours. Let us therefore first derive some results by hand. For a lower bound on
the right-hand tail of $Z$, no heavy machinery is necessary. Indeed let
$i_{0}=\arg\max_{i}K_{ii}$ and $\mu=E\big(\max N_{i}\big)>0$. Then, for
$z>0$,
\begin{equation}
P\big(Z\geqslant z\big)\geqslant P\big(N_{i_{0}}\geqslant\mu+z\big)\geqslant
\frac{1}{\sqrt{2\pi}}~\frac{(\mu+z)^{2}}{\sigma_{\max}^{2}+(\mu+z)^{2}%
}~e^{-\frac{(\mu+z)^{2}}{2\,\sigma_{\max}^{2}}}. \label{trivlb}%
\end{equation}
A nearly identical argument leads to the following upper bound on the
left-hand tail of $Z$: for $z>0$,
\begin{equation}
P(Z\leqslant-z)\leqslant\frac{\min_{i}\sqrt{K_{ii}}}{\sqrt{2\pi}(z-\mu)}%
\exp\left(  -\frac{(z-\mu)^{2}}{2\min_{i}K_{ii}}\right)  . \label{trivub}%
\end{equation}
This improves Borell's inequality (\ref{borell-b}) asymptotically.

By using the techniques in our article, the density estimates in
(\ref{rho-intro}) allow us to obtain a new lower bound result on $Z$'s left
hand tail, and to improve the classical right-hand tail result of
(\ref{vitale-b}). We have for the right-hand tail
\begin{equation}
\frac{E|Z|\,\sigma_{\min}^{2}}{2\,\sigma_{\max}^{2}}~\frac{z}{\sigma_{\min
}^{2}+z^{2}}~\mathrm{exp}\left(  -\frac{z^{2}}{2\,\sigma_{\min}^{2}}\right)
\leqslant P\big(Z\geqslant z\big)\leqslant\frac{E|Z|\,\sigma_{\max}^{2}%
}{2\,\sigma_{\min}^{2}}~\frac{1}{z}~\mathrm{exp}\left(  -\frac{z^{2}%
}{2\,\sigma_{\max}^{2}}\right)  \label{improve-intro}%
\end{equation}
if $z>0$, and one notes that the above right-hand side goes (slightly) faster
to zero than (\ref{vitale-borell-b}), because of the presence of the factor
$z^{-1}$; yet the lower bound is less sharp than (\ref{trivlb}) for large $z$.
The first and last expressions in (\ref{improve-intro}) are also lower and
upper bounds for the left-hand tail $P\big(Z\leqslant-z\big)$. To the best of
our knowledge, the lower bound is new; the upper bound is less sharp than
(\ref{trivub}) for large $z$.\vspace{0.1in}

Let us now cite some works which are related to ours, insofar as some of the
preoccupations and techniques are similar. I{n \cite{houdreprivault},
Houdr\'{e} and Privault prove concentration inequalities for functionals of
Wiener and Poisson spaces: they have discovered almost-sure conditions on
expressions involving Malliavin derivatives which guarantee upper bounds on
the tails of their functionals. This is similar to the upper bound portion of
our work in Section 4, and closer yet to the first-chaos portion of the work
in \cite{VVJFA}; they do not, however, address lower bound issues, nor do they
have any claims regarding densities. 
%%Je préfère enlever ce passage, car la fin est un peu négative vis-à-vis du 
%%travail d'Houdré et Privault. Et comme j'aime beaucoup Houdré, j'aimerais autant
%%ne pas me fâcher avec lui! (s'il lit ça, je suis en effet sûr qu'il va m'écrire un message
%%pour me dire que ce n'est pas vrai car ceci cela...)
%%
%%Their work unifies the treatment of
%%Malliavin calculus in the Wiener and Poisson spaces, and goes quite far in the
%%analysis of the Poisson case; but when restricted to the Wiener space, which
%%is our preoccupation, \cite{houdreprivault} appears not to exceed the now
%%classical analysis in \cite{pisier}.}

Decreusefond and Nualart \cite{decreununu} obtain, by means of the Malliavin
calculus, estimates for the Laplace transform of the hitting times of any
general Gaussian process; they define a monotonicity condition on the
covariance function of such a process under which this Laplace transform is
bounded above by that of standard Brownian motion; similarly to how we derive
upper tail estimates of Gaussian type from our analysis, they derive the
finiteness of some moments by comparison to the Brownian case. However, as in
\cite{houdreprivault}, {reference} \cite{decreununu} does not address issues
of densities or of lower bounds.

General lower bound results on densities are few and far between. The case of
uniformly elliptic diffusions was treated in a series of papers by Kusuoka and
Stroock: see \cite{KS}. This was generalized by Kohatsu-Higa \cite{K} in
Wiener space via the concept of uniformly elliptic random variables; these
random variables proved to be well-adapted to studying diffusion equations.
%%but were otherwise relatively unwieldy. 
E. Nualart \cite{Ncras} showed that
fractional exponential moments for a divergence-integral quantity known to be
useful for bounding densities from above (see formula (\ref{rho}) below), can
also be useful for deriving a scale of exponential lower bounds on densities;
the scale includes Gaussian lower bounds. However, in all these works, the
applications are largely restricted to diffusions.\vspace{0.1in}

We now introduce our general setting which will allow to prove
(\ref{rho-intro})-(\ref{ffbbmm}) and several other results. We consider a
centered isonormal Gaussian process $X=\{X(h):h\in\EuFrak H\}$ defined on a
real separable Hilbert space $\EuFrak H$. This just means that $X$ is a
collection of centered and jointly Gaussian random variables indexed by the
elements of $\EuFrak H$, defined on some probability space $(\Omega
,\mathscr{F},P)$ and such that, for every $h,g\in\EuFrak H$,
\[
E\big(X(h)X(g)\big)=\langle h,g\rangle_{\EuFrak H}.
\]
As usual in Malliavin calculus, we use the following notation (see Section
\ref{preliminaries} for precise definitions):

\begin{itemize}
\item $L^{2}(\Omega,\mathscr{F},P)$ is the space of square-integrable
functionals of $X$. This means in particular that $\mathscr{F}$ is the
$\sigma$-field generated by $X$;

\item ${\mathbb{D}}^{1,2}$ is the domain of the Malliavin derivative operator
$D$ with respect to $X$. {Roughly speaking, it} is the subset of random
variables in $L^{2}(\Omega,\mathscr{F},P)$ whose Malliavin derivative is also
in $L^{2}(\Omega,\mathscr{F},P)$;

\item $\mathrm{Dom}\delta$ is the domain of the divergence operator $\delta$.
This operator will {really }only play a marginal role in our study; it is
{simply }used in order to simplify some proof arguments, and for comparison purposes.
\end{itemize}

From now on, $Z$ will always denote a random variable of $\mathbb{D}^{1,2}$
with \textit{zero mean}. Recall that its derivative $DZ$ is a random element
with values in $\EuFrak H$. The following result on the density of a random
variable is a well-known fact of the Malliavin calculus: if $DZ/\Vert
DZ\Vert_{\EuFrak H}^{2}$ belongs to $\mathrm{Dom}\delta$, then $Z$ has a
continuous and bounded density $\rho$ given, for all $z\in\mathbb{R}$, by
\begin{equation}
\rho(z)=E\left[  \mathbf{1}_{(z,+\infty]}(Z)\,\delta\left(  \frac{DZ}{\Vert
DZ\Vert_{\EuFrak H}^{2}}\right)  \right]  . \label{rho}%
\end{equation}
From this expression, it is sometimes possible to deduce \textit{upper}
{bounds} for $\rho$. Several examples are detailed in Section 2.1.1 of
Nualart's book \cite{nualartbook}. Note the following two points, however: (a)
it is not clear whether it is at all possible to prove (\ref{rho-intro}) by
using (\ref{rho}); (b) more generally it appears to be just as difficult to
deduce any \textit{lower}-bound\ relations on the density $\rho$ of any random
variable via (\ref{rho}).\vspace{0.1in}

Herein we prove a new general formula for $\rho$, from which we easily deduce
(\ref{rho-intro}) for instance. For $Z$ a mean-zero r.v. in $\mathbb{D}^{1,2}%
$, define the function $g:\mathbb{R\to R}$ 
%%J'ai rajouté:
almost everywhere 
by
\begin{equation}
g(z)  =g_{Z}(z):=E\left(  \left.  \langle
DZ,-DL^{-1}Z\rangle_{\EuFrak H}\right\vert Z=z\right)  . \label{gee}%
\end{equation}
The $L$ appearing here is the so-called {generator of the Ornstein-Uhlenbeck
semigroup}, defined in the next section. We drop the subscript $_{Z}$ from
$g_{Z}$ in this article, since each example herein refers to only one r.v. $Z$
at a time. By \cite[Proposition 3.9]{NP}, $g$ is non-negative on the support
of $Z$. Under some general conditions on $Z$ (see Theorem \ref{key-thm} for a
precise statement), the density $\rho$ of $Z$ is given by the following new
formula, for any $z$ in $Z$'s support:
%%J'ai un peu changé la manière d'écrire
\begin{equation}
P\left(  Z\in dz\right)  =\rho(z)dz=\frac{E|Z|}{2g(z)
}~\mathrm{exp}\left(  -\int_{0}^{z}\frac{x\,dx}{g(x)}\right)dz.
\label{rho2}%
\end{equation}
\vspace{0.1in}

The key point in our approach is that it is possible, in many cases, to
estimate the quantity $g(z)  $ in (\ref{gee}) rather precisely.
In particular, we will make systematic use of the following consequence of the
%%J'ai enlevé "so called"
Mehler formula (see Remark 3.6 in \cite{NP}),
%%J'ai rajouté "see Remark 3.6 in \cite{NP} ci dessus
also proved herein (Proposition \ref{computation}):%
\[
g(z)  =\int_{0}^{\infty}e^{-u}\,\mathbf{E}\big(\langle\Phi
_{Z}(X),\Phi_{Z}(e^{-u}X+\sqrt{1-e^{-2u}}X^{\prime})\rangle_{\EuFrak
H}|Z=z\big)du.
\]
In this formula, the mapping $\Phi_{Z}:\mathbb{R}^{\EuFrak H}\rightarrow
\EuFrak H$ is defined $P\circ X^{-1}$-almost surely through the identity
$DZ=\Phi_{Z}(X)$, while $X^{\prime}$, which stands for an independent copy of
$X$, is such that $X$ and $X^{\prime}$ are defined on the product probability
space $(\Omega\times\Omega^{\prime},\mathscr{F}\otimes\mathscr{F}^{\prime
},P\times P^{\prime})$; $\mathbf{E}$ denotes the mathematical expectation with
respect to $P\times P^{\prime}$. 
%%J'ai changé un peu la phrase ci-dessus
This formula for $g$ then allows, in many cases,
to obtain via (\ref{rho2}) a lower and an upper bound on $\rho$ simultaneously.
%%Fin changement
We refer the reader to Corollary \ref{main-thm} and the examples in Section
\ref{density}, and in particular to {the second and fourth examples, which are
the proofs of Proposition \ref{rhointroprop}} and Proposition \ref{pr}
respectively. At this stage, let us note however that it is not possible to
obtain only a lower bound, or only an upper bound, using formula (\ref{rho2}).
Indeed, one can see that one needs to control $g$ simultaneously from above
and below to get the technique to work.\bigskip

In the second main part of the paper (Section \ref{sect-concentration}), we
explain what can be done when one only knows how to bound $g$ from one
direction or the other, but not both simultaneously. Note that one is
precisely in this situation when one seeks to prove the inequalities
(\ref{vitale-b}) and (\ref{borell-b}). These will be a simple consequence of a
more general upper bound result (Theorem \ref{concentration}) in Section
\ref{sect-concentration}.

As another application of Theorem \ref{concentration}, the following result
concerns a functional of fractional Brownian motion.

\begin{prop}
\label{fbm} Let $B=\{B_{t},\,t\in\lbrack0,T]\}$ be a fractional Brownian motion with
Hurst index $H\in(0,1)$. Then, {{denoting $c_{H}=H+1/2$, we have}}, for any
$z>0$:%
\[
P\left(  \int_{0}^{T}B_{u}^{2}du\geqslant z+T^{2H+1}/(2c_{H})\right)
\leqslant\mathrm{exp}\left(  {{-\frac{c_{H}^{2}~z^{2}}{2c_{H}T^{2H+1}%
z+T^{4H+2}}}}\right)  .
\]

\end{prop}

Of course, the interest of this result lies in the fact that the exact
distribution of $\int_{0}^{T}B_{u}^{2}du$ is still an open problem {when
$H\neq1/2$}. With respect to the classical result by Borell \cite{borell2}
(which would give a bound like $\mathrm{exp}(-Cz)$), observe here that, {as in
Chatterjee \cite{Chatterbox},} we get a kind of \textquotedblleft
continuous\textquotedblright\ transition from Gaussian to exponential tails.
The behavior for large $z$ is always of exponential type. At the end of this
article, we take up the issue of finding a lower bound which might be
commensurate with the upper bound above; our Malliavin calculus techniques
fail here, but we are still able to derive an interesting result by hand, see
(\ref{lbfBm}).\vspace{0.1in}

Section \ref{sect-concentration} also contains a lower bound result, Theorem
\ref{thmlbnostein}, again based on the quantity $\langle DZ,-DL^{-1}%
Z\rangle_{\EuFrak H}$ via the function $g$ in (\ref{gee}). This quantity was
introduced recently in \cite{NP} for the purpose of {{using Stein's method in
order to show that the standard deviation of $\langle DZ,-DL^{-1}%
Z\rangle_{\EuFrak H}$ provides an error bound of the normal approximation of
$Z$, see also Remark \ref{ivanjoerem} below. Here, in Theorem
\ref{thmlbnostein} and in Theorem \ref{concentration} as a special case
($\alpha=0$ therein), $g(Z)  =E(\langle DZ,-DL^{-1}%
Z\rangle_{\EuFrak H}|Z)$ can be instead assumed to be bounded either above {\it or}
below {\textit{almost surely }}by a constant; this constant's role is to be a
measure of the variance of $Z$, and more specifically to ensure that the tail
of }}$Z$ is bounded {either above or below by a normal tail with that constant
as its variance.} Our Section \ref{sect-concentration} {{can thus be thought}
as a way to extend the phenomena described in \cite{NP} when comparison with
the normal distribution can only be expected to go one way. Theorem
\ref{thmlbnostein} shows that we may have no control over how heavy the tail
of $Z$ may be (beyond the existence of a second moment), but {the condition}
$g(Z)  \geqslant{{\sigma^{2}>0}}$ essentially guarantees that it
has to be no less heavy than a Gaussian tail with variance }$\sigma^{2}$%
{.}\vspace{0.1in}

We finish this description of our results by stressing again that, whether in
Sections \ref{density} or \ref{sect-concentration}, we present many examples
where the quantities $\langle DZ,-DL^{-1}Z\rangle_{\EuFrak H}$ and
$\langle\Phi_{Z}(X),\Phi_{Z}(e^{-u}X+\sqrt{1-e^{-2u}}X^{\prime})\rangle
_{\EuFrak H}$ are computed and estimated easily, by hand and/or via
Proposition \ref{computation}. The advantage over formulas such as
(\ref{rho}), which involve the unwieldy divergence operator $\delta$, should
be clear.\vspace{0.1in}

The rest of the paper is organized as follows. In Section 2, we recall the
notions of Malliavin calculus that we need in order to perform our proofs. In
Section 3, we state and discuss our density estimates. Section 4 deals with
concentration inequalities, i.e. tail estimates.

\section{Some elements of Malliavin calculus}

\label{preliminaries} We follow Nualart's book \cite{nualartbook}. As stated
in the introduction, we denote by $X$ a centered isonormal Gaussian process over a
real separable Hilbert space $\EuFrak H$. Let $\mathscr{F}$ be the $\sigma
$-field generated by $X$. It is well-known that any random variable $Z$
belonging to $L^{2}(\Omega,\mathscr{F},P)$ admits the following chaos
expansion:
\begin{equation}
Z=\sum_{m=0}^{\infty}I_{m}(f_{m}), \label{ChaosExpansion}%
\end{equation}
where $I_{0}(f_{0})=E(Z)$, the series converges in $L^{2}(\Omega)$ and the
kernels $f_{m}\in\EuFrak H^{\odot m}$, $m\geqslant1$, are uniquely determined
by $Z$. In the particular case where $\EuFrak H=L^{2}(A,\mathscr{A},\mu)$, for
$(A,\mathscr{A})$ a measurable space and $\mu$ a $\sigma$-finite and
non-atomic measure, one has that $\EuFrak H^{\odot m}=L_{s}^{2}(A^{m}%
,\mathscr{A}^{\otimes m},\mu^{\otimes m})$ is the space of symmetric and
square integrable functions on $A^{m}$ and, for every $f\in\EuFrak H^{\odot
m}$, $I_{m}(f)$ coincides with the multiple Wiener-It\^{o} integral of order
$m$ of $f$ with respect to $X$. For every $m\geqslant0$, we write $J_{m}$ to
indicate the orthogonal projection operator on the $m$th Wiener chaos
associated with $X$. That is, if $Z\in L^{2}(\Omega,\mathscr{F},P)$ is as in
(\ref{ChaosExpansion}), then $J_{m}F=I_{m}(f_{m})$ for every $m\geqslant0$.

Let $\mathscr{S}$ be the set of all smooth cylindrical random variables of the
form
\[
Z=g\big(X(\phi_{1}),\ldots,X(\phi_{n})\big)
\]
where $n\geqslant1$, $g:\mathbb{R}^{n}\rightarrow\mathbb{R}$ is a smooth
function with compact support and $\phi_{i}\in\EuFrak H$. The Malliavin
derivative of $Z$ with respect to $X$ is the element of $L^{2}(\Omega
,\EuFrak
H)$ defined as
\[
DZ=\sum_{i=1}^{n}\frac{\partial g}{\partial x_{i}}\big(X(\phi_{1}%
),\ldots,X(\phi_{n})\big)\phi_{i}.
\]
In particular, $DX(h)=h$ for every $h\in\EuFrak H$. By iteration, one can
define the $m$th derivative $D^{m}Z$ (which is an element of $L^{2}%
(\Omega,\EuFrak H^{\odot m})$) for every $m\geqslant2$. As usual, for
$m\geqslant1$, ${\mathbb{D}}^{m,2}$ denotes the closure of $\mathscr{S}$ with
respect to the norm $\Vert\cdot\Vert_{m,2}$, defined by the relation
\[
\Vert Z\Vert_{m,2}^{2}=E(Z^{2})+\sum_{i=1}^{m}E\big(\Vert D^{i}Z\Vert
_{\EuFrak
H^{\otimes i}}^{2}\big).
\]
Note that a random variable $Z$ as in (\ref{ChaosExpansion}) is in
${\mathbb{D}}^{1,2}$ if and only if
\[
\sum_{m=1}^{\infty}m\,m!\,\Vert f_{m}\Vert_{\EuFrak H^{\otimes m}}^{2}%
<\infty,
\]
and, in this case, $E\big(\Vert DZ\Vert_{\EuFrak H}^{2}\big)=\sum
_{m\geqslant1}m\,m!\,\Vert f_{m}\Vert_{\EuFrak H^{\otimes m}}^{2}$. If
$\EuFrak H=L^{2}(A,\mathscr{A},\mu)$ (with $\mu$ non-atomic), then the
derivative of a random variable $Z$ as in (\ref{ChaosExpansion}) can be
identified with the element of $L^{2}(A\times\Omega)$ given by
\[
D_{a}Z=\sum_{m=1}^{\infty}mI_{m-1}\big(f_{m}(\cdot,a)\big),\quad a\in A.
\]
The Malliavin derivative $D$ satisfies the following chain rule. If
$\varphi:\mathbb{R}^{n}\rightarrow\mathbb{R}$ is of class $\mathscr{C}^{1}$
with bounded derivatives, and if $\{Z_{i}\}_{i=1,\ldots,n}$ is a vector of
elements of ${\mathbb{D}}^{1,2}$, then $\varphi(Z_{1},\ldots,Z_{n}%
)\in{\mathbb{D}}^{1,2}$ and
\begin{equation}
D\,\varphi(Z_{1},\ldots,Z_{n})=\sum_{i=1}^{n}\frac{\partial\varphi}{\partial
x_{i}}(Z_{1},\ldots,Z_{n})DZ_{i}. \label{chainrule-lipschitz}%
\end{equation}
Formula (\ref{chainrule-lipschitz}) still holds when $\varphi$ is only
Lipshitz but the law of $(Z_{1},\ldots,Z_{n})$ has a density with respect to
the Lebesgue measure on $\mathbb{R}^{n}$ (see e.g. Proposition 1.2.3 in
\cite{nualartbook}).

We denote by $\delta$ the adjoint of the operator $D$, also called the
divergence operator. A random element $u\in L^{2}(\Omega,\EuFrak H)$ belongs
to the domain of $\delta$, denoted by $\mathrm{Dom}\delta$, if and only if it
satisfies
\[
\big|E\langle DZ,u\rangle_{\EuFrak H}\big|\leqslant c_{u}\,E(Z^{2})^{1/2}%
\quad\mbox{for any }Z\in{\mathscr{S}},
\]
where $c_{u}$ is a constant depending only on $u$. If $u\in\mathrm{Dom}\delta
$, then the random variable $\delta(u)$ is uniquely defined by the duality
relationship
\begin{equation}
E(Z\delta(u))=E\langle DZ,u\rangle_{\EuFrak H}, \label{ipp}%
\end{equation}
which holds for every $Z\in{\mathbb{D}}^{1,2}$.

The operator $L$ is defined through the projection operators as $L=\sum
_{m=0}^{\infty}-mJ_{m},$ and is called the generator of the Ornstein-Uhlenbeck
semigroup. It satisfies the following crucial property. A random variable $Z$
is an element of $\mathrm{Dom}$$L\,\,(={\mathbb{D}}^{2,2})$ if and only if
$Z\in\mathrm{Dom}\delta D$ (i.e. $Z\in{\mathbb{D}}^{1,2}$ and $DZ\in
\mathrm{Dom}\delta$), and in this case:
\begin{equation}
\delta DZ=-LZ. \label{deltaD=L}%
\end{equation}
We also define the operator $L^{-1}$, which is the inverse of $L$, as follows.
For every $Z\in L^{2}(\Omega,\mathscr{F},P)$, we set $L^{-1}Z$ $=$
$\sum_{m\geqslant1}-\frac{1}{m}J_{m}(Z)$. Note that $L^{-1}$ is an operator
with values in ${\mathbb{D}}^{2,2}$, and that $LL^{-1}Z=Z-E(Z)$ for any $Z\in
L^{2}(\Omega,\mathscr{F},P)$, so that $L^{-1}$ does act as $L$'s inverse for
centered r.v.'s.

The family $(T_{u},\,u\geqslant0)$ of operators is defined as $T_{u}%
=\sum_{m=0}^{\infty}e^{-mu}J_{m}$, and is called the Orstein-Uhlenbeck
semigroup. Assume that the process $X^{\prime}$, which stands for an
independent copy of $X$, is such that $X$ and $X^{\prime}$ are defined on the
product probability space $(\Omega\times\Omega^{\prime},\mathscr{F}\otimes
\mathscr{F}^{\prime},P\times P^{\prime})$. Given a random variable
$Z\in\mathbb{D}^{1,2}$, we can write $DZ=\Phi_{Z}(X)$, where $\Phi_{Z}$ is a
measurable mapping from $\mathbb{R}^{\EuFrak H}$ to $\EuFrak H$, determined
$P\circ X^{-1}$-almost surely. Then, for any $u\geqslant0$, we have the
so-called Mehler formula:
\begin{equation}
T_{u}(DZ)=E^{\prime}\big(\Phi_{Z}(e^{-u}X+\sqrt{1-e^{-2u}}X^{\prime})\big),
\label{mehler}%
\end{equation}
where $E^{\prime}$ denotes the mathematical expectation with respect to the
probability $P^{\prime}$.

\section{Density estimates}

\label{density}

For $Z\in\mathbb{D}^{1,2}$ with zero mean, recall the function $g$ introduced in the
introduction in (\ref{gee}):%
\[
g(z)  =E(\langle DZ,-DL^{-1}Z\rangle_{\EuFrak H}|Z=z).
\]
It is useful to keep in mind throughout this paper that, by \cite[Proposition
3.9]{NP}, $g(z)  \geqslant0$ on the support of $Z$. In this
section, we further assume that $g$ is bounded away from $0$.

\subsection{General formulae and estimates}

We begin with the following theorem, which will be key in the sequel.

\begin{thm}
\label{key-thm} Let $Z\in\mathbb{D}^{1,2}$ with zero mean, and $g$ as above.
Assume that there exists $\sigma_{\min}>0$ such that
\begin{equation}
g(Z)  \geqslant\sigma_{\min}^{2}\quad\mbox{almost surely.}
\label{sigmamin}%
\end{equation}
Then $Z$ has a density $\rho$, its support is $\mathbb{R}$ and we have, almost
everywhere:
\begin{equation}
\rho(z)=\frac{E\left\vert Z\right\vert }{2g(z)  }\,\,\mathrm{exp}%
\left(  -\int_{0}^{z}\frac{x\,dx}{g(x)  }\right)  .
\label{ohquellebelleformulebis}%
\end{equation}

\end{thm}

\textit{Proof}. We split the proof into several steps.\newline

\textit{Step 1: An integration by parts formula}. For any $f:\mathbb{R}%
\rightarrow\mathbb{R}$ of class $\mathscr{C}^{1}$ with bounded derivative, we
have
\begin{align}
E\big(Zf(Z)\big)  &  =E\big(LL^{-1}Zf(Z)\big)=E\big(\delta D(-L^{-1}%
Z)f(Z)\big)\quad\mbox{by (\ref{deltaD=L})}\nonumber\\
&  =E\big(\langle Df(Z),-DL^{-1}Z\rangle_{\EuFrak H}\big)\quad\mbox{by
(\ref{ipp})}\nonumber\\
&  =E\big(f^{\prime}(Z)\langle DZ,-DL^{-1}Z\rangle_{\EuFrak H}\big)\quad
\mbox{by (\ref{chainrule-lipschitz})}. \label{eheh}%
\end{align}

\vskip0.3cm

\textit{Step 2: Existence of the density}. Fix $a<b$ in $\mathbb{R}$. For any
$\varepsilon>0$, consider a $\mathscr{C}^{\infty}$-function $\varphi
_{\varepsilon}:\mathbb{R}\rightarrow\lbrack0,1]$ such that $\varphi
_{\varepsilon}(z)=1$ if $z\in\lbrack a,b]$ and $\varphi_{\varepsilon}(z)=0$ if
$z<a-\varepsilon$ or $z>b+\varepsilon$. We set $\psi_{\varepsilon}%
(z)=\int_{-\infty}^{z}\varphi_{\varepsilon}(y)dy$ for any $z\in\mathbb{R}$.
Then, we can write
\begin{align*}
P(a\leqslant Z\leqslant b)  &  =E\big(\mathbf{1}_{[a,b]}(Z)\big)\\
&  \leqslant\sigma_{\min}^{-2}\,E\big(\mathbf{1}_{[a,b]}(Z)E(\langle
DZ,-DL^{-1}Z\rangle_{\EuFrak H}|Z)\big)\quad\mbox{by assumption
(\ref{sigmamin})}\\
&  =\sigma_{\min}^{-2}\,E\big(\mathbf{1}_{[a,b]}(Z)\langle DZ,-DL^{-1}%
Z\rangle_{\EuFrak H}\big)\\
&  =\sigma_{\min}^{-2}\,E\big(\liminf_{\varepsilon\rightarrow0}\varphi
_{\varepsilon}(Z)\langle DZ,-DL^{-1}Z\rangle_{\EuFrak H}\big)\\
&  \leqslant\sigma_{\min}^{-2}\,\liminf_{\varepsilon\rightarrow0}%
E\big(\varphi_{\varepsilon}(Z)\langle DZ,-DL^{-1}Z\rangle_{\EuFrak H}%
\big)\quad\mbox{by
Fatou's inequality}\\
&  =\sigma_{\min}^{-2}\,\liminf_{\varepsilon\rightarrow0}E\big(\psi
_{\varepsilon}(Z)Z\big)\quad\mbox{by (\ref{eheh})}\\
&  =\sigma_{\min}^{-2}\,E\left(  Z\int_{-\infty}^{Z}\mathbf{1}_{[a,b]}%
(u)du\right)  \quad\mbox{by bounded
convergence}\\
&  =\sigma_{\min}^{-2}\int_{a}^{b}E\big(Z\mathbf{1}_{[u,+\infty)}%
(Z)\big)du\leqslant(b-a)\times\sigma_{\min}^{-2}\,E|Z|.
\end{align*}
This implies the absolute continuity of $Z$, that is the existence of $\rho$.

\vskip0.3cm

\textit{Step 3: A key formula}. Let $f:\mathbb{R}\rightarrow\mathbb{R}$ be a
continuous function with compact support, and $F$ denote any antiderivative of
$f$. Note that $F$ is bounded. We have
\begin{align}
E\big(f(Z)\langle DZ,-DL^{-1}Z\rangle_{\EuFrak H}\big)  &
=E\big(F(Z)Z\big)\quad\mbox{by (\ref{eheh})}\nonumber\\
&  =\int_{\mathbb{R}}F(z)\,z\,\rho(z)dz\nonumber\\
&  \underset{(\ast)}{=}\int_{\mathbb{R}}f(z)\left(  \int_{z}^{\infty}%
y\rho(y)dy\right)  dz\nonumber\\
&  =E\left(  f(Z)\frac{\int_{Z}^{\infty}y\rho(y)dy}{\rho(Z)}\right)
.\nonumber
\end{align}
Equality (*) was obtained by integrating by parts, after observing that
\[
\int_{z}^{\infty}y\rho(y)dy\longrightarrow0\,\mbox{ as $|z|\to\infty$}
\]
(for $z\rightarrow+\infty$, this is because $Z\in L^{1}(\Omega)$; for
$z\rightarrow-\infty$, this is because $Z$ has mean zero). Therefore, we have
shown
\begin{equation}
g(Z)  =E(\langle DZ,-DL^{-1}Z\rangle_{\EuFrak H}|Z)=\frac
{\int_{Z}^{\infty}y\rho(y)dy}{\rho(Z)}\quad\mbox{almost surely}.
\label{ohlabelleformule}%
\end{equation}

\textit{Step 4: The support of $\rho$}. Since $Z\in\mathbb{D}^{1,2}$, it is
known (see e.g. \cite[Proposition 2.1.7]{nualartbook}) that $\mathrm{Supp}%
\rho=[\alpha,\beta]$ with $-\infty\leqslant\alpha<\beta\leqslant+\infty$.
Since $Z$ has zero mean, note that $\alpha<0$ and $\beta>0$ necessarily.
Identity (\ref{ohlabelleformule}) yields
\begin{equation}
\int_{z}^{\infty}y\rho\left(  y\right)  dy\geqslant\sigma_{\min}^{2}%
\,\rho\left(  z\right)  \quad\mbox{for almost all $z\in(\alpha,\beta)$}.
\label{fgr}%
\end{equation}
For every $z\in(\alpha,\beta)$, define $\varphi\left(  z\right)  :=\int
_{z}^{\infty}y\rho\left(  y\right)  dy.$ This function is differentiable
almost everywhere on $(\alpha,\beta)$, and its derivative is $-z\rho\left(
z\right)  $. In particular, since $\varphi(\alpha)=\varphi(\beta)=0$, we have
that $\varphi(z)>0$ for all $z\in(\alpha,\beta)$. On the other hand,
when multiplied by $z\in\lbrack{0,\beta)}$, the inequality (\ref{fgr}) gives
$\frac{\varphi^{\prime}\left(  z\right)  }{\varphi\left(  z\right)  }%
\geqslant-\frac{z}{\sigma_{\min}^{2}}$. Integrating this relation over the
interval $[0,z]$ yields $\log\varphi\left(  z\right)  -\log\varphi\left(
0\right)  \geqslant-\frac{z^{2}}{2\,\sigma_{\min}^{2}}$, i.e., since
$0=E(Z)=E(Z_{+})-E(Z_{-})$ so that $E|Z|=E(Z_{+})+E(Z_{-})=2E(Z_{+}%
)=2\varphi(0)$, we have
\begin{equation}
\varphi\left(  z\right)  =\int_{z}^{\infty}y\rho\left(  y\right)
dy\geqslant\frac{1}{2}E|Z|e^{-\frac{z^{2}}{2\,\sigma_{\min}^{2}}}. \label{Flb}%
\end{equation}
Similarly, when multiplied by $z\in(\alpha,0]$, inequality (\ref{fgr}) gives
$\frac{\varphi^{\prime}\left(  z\right)  }{\varphi\left(  z\right)  }%
\leqslant-\frac{z}{\sigma_{\min}^{2}}.$ Integrating this relation over the
interval $[z,0]$ yields $\log\varphi\left(  0\right)  -\log\varphi\left(
z\right)  \leqslant\frac{z^{2}}{2\,\sigma_{\min}^{2}}$, i.e. (\ref{Flb}) still
holds for $z\in(\alpha,0]$. Now, let us prove that $\beta=+\infty$. If this
were not the case, by definition, we would have $\varphi\left(  \beta\right)
=0$; on the other hand, by letting $z$ tend to $\beta$ in the above
inequality, because $\varphi$ is continuous, we would have $\varphi\left(
\beta\right)  \geqslant\frac{1}{2}E|Z|e^{-\frac{\beta^{2}}{2\sigma_{\min}^{2}%
}}>0$, which contradicts $\beta<+\infty$. The proof of $\alpha=-\infty$ is
similar. In conclusion, we have shown that $\mathrm{supp}\rho=\mathbb{R}$.

\vskip0.3cm

\textit{Step 5: Proof of (\ref{ohquellebelleformulebis})}. Let $\varphi
:\mathbb{R}\rightarrow\mathbb{R}$ be still defined by $\varphi(z)=\int
_{z}^{\infty}y\rho(y)dy$. On one hand, we have $\varphi^{\prime}(z)=-z\rho(z)$
for almost all $z\in\mathbb{R}$. On the other hand, by (\ref{ohlabelleformule}%
), we have, for almost all $z\in\mathbb{R}$,
\begin{equation}
\varphi(z)=\rho(z)g(z)  . \label{stop1}%
\end{equation}
By putting these two facts together, we get the following ordinary
differential equation satisfied by $\varphi$:
\[
\frac{\varphi^{\prime}(z)}{\varphi(z)}=-\frac{z}{g(z)  }%
\quad\mbox{for almost all $z\in\R$.}
\]
Integrating this relation over the interval $[0,z]$ yields
\[
\log\varphi(z)=\log\varphi(0)-\int_{0}^{z}\frac{x\,dx}{g(x)  }.
\]
Taking the exponential and using the fact that $\varphi(0)=\frac{1}{2}\,E|Z|$,
we get
\[
\varphi(z)=\frac{1}{2}\,E|Z|\,\mathrm{exp}\left(  -\int_{0}^{z}\frac
{x\,dx}{g(x)  }\right)  .
\]
Finally, the desired conclusion comes from (\ref{stop1}). \vspace{-0.6cm}
\begin{flushright}
\mbox{$\Box$}
\end{flushright}\noindent

\begin{rem}
\label{ivanjoerem} {\textrm{The \textquotedblleft integration by parts
formula\textquotedblright\ (\ref{eheh}) was proved and used for the first time
by Nourdin and Peccati in \cite{NP}, in order to perform error bounds in the
normal approximation of $Z$. Specifically, \cite{NP} shows, by combining
Stein's method with (\ref{eheh}), that
\begin{equation}
\sup_{z\in\mathbb{R}}\big|P(Z\leqslant z)-P(N\leqslant z)\big|\leqslant
\frac{\sqrt{{\rm Var}\big(E(\langle DZ,-DL^{-1}Z\rangle_{\EuFrak H}|Z)\big)}%
}{{\rm Var}(Z)},\label{ivanjoe}%
\end{equation}
where $N\sim\mathscr{N}(0,{\rm Var}Z)$. In reality, the inequality stated in
\cite{NP} is with ${\rm Var}\big(\langle DZ,-DL^{-1}Z\rangle_{\EuFrak H}\big)$
instead of ${\rm Var}\big(E(\langle DZ,-DL^{-1}Z\rangle_{\EuFrak H}|Z)\big)$ on the
right-hand side; but the same proof allows to write this slight improvement;
it was not stated or used in \cite{NP} because it did not improve the
applications therein.}}
\end{rem}

Using Theorem \ref{key-thm}, we can deduce the following interesting criterion
for normality, which one will compare with (\ref{ivanjoe}).

\begin{cor}
\label{levy} Let $Z\in\mathbb{D}^{1,2}$; let $g(Z)  =E(\langle
DZ,-DL^{-1}Z\rangle_{\EuFrak H}|Z)$. Then $Z$ is Gaussian if and only if
$\mathrm{Var}(g(Z)  )=0.$
\end{cor}

\textit{Proof}: We can assume without loss of generality that $Z$ is centered.
By (\ref{eheh}) (choose $f(z)=z$), we have
\[
E(\langle DZ,-DL^{-1}Z\rangle_{\EuFrak H})=E(Z^{2})=\mathrm{Var}Z.
\]
Therefore, the condition $\mathrm{Var}(g(Z)  )=0$ is equivalent
to
\[
g(Z)  =\mathrm{Var}Z\quad\mbox{almost
surely.}
\]
Let $Z\sim\mathscr{N}(0,\sigma^{2})$. Using (\ref{ohlabelleformule}), we
immediately check that $g(Z)  =\sigma^{2}$ almost surely.
Conversely, if $g(Z)  =\sigma^{2}$ almost surely, then Theorem
\ref{key-thm} implies that $Z$ has a density $\rho$ given by $\rho
(z)=\frac{E|Z|}{2\sigma^{2}}e^{-\frac{z^{2}}{2\,\sigma^{2}}}$ for almost all
$z\in\mathbb{R}$, from which we immediately deduce that $Z\sim
\mathscr{N}(0,\sigma^{2})$. \vspace{-0.6cm} \begin{flushright}
\mbox{$\Box$}
\end{flushright}\noindent

Observe that if $Z\sim\mathscr{N}(0,\sigma^{2})$, then $E|Z|=\sqrt{2/\pi
}\,\sigma$, so that the formula (\ref{ohquellebelleformulebis}) for $\rho$
agrees, of course, with the usual one in this case.\newline

Depending on the situation, $g(Z)  $ may be computable or may be
estimated by hand. We cite the next corollary for situations where this is the
case. However, with the exception of this corollary, the remainder of this
section, starting with Proposition \ref{computation}, provides a systematic
computational technique to deal with $g(Z)  $.

\begin{cor}
\label{general} Let $Z\in\mathbb{D}^{1,2}$ with zero mean and 
$g(Z)  :=E(\langle DZ,-DL^{-1}Z\rangle_{\EuFrak H}|Z)$. If there exists
$\sigma_{\min},\sigma_{\max}>0$ such that
\[
\sigma_{\min}^{2}\leqslant g(Z)  \leqslant\sigma_{\max}^{2}%
\quad\mbox{almost surely},
\]
then $Z$ has a density $\rho$ satisfying, for almost all $z\in\mathbb{R}$%
\[
\frac{E|Z|}{2\,\sigma_{\min}^{2}}\,\mathrm{exp}\left(  -\frac{z^{2}}%
{2\sigma_{\max}^{2}}\right)  \leqslant\rho(z)\leqslant\frac{E|Z|}%
{2\,\sigma_{\max}^{2}}\,\mathrm{exp}\left(  -\frac{z^{2}}{2\sigma_{\min}^{2}%
}\right)  .
\]

\end{cor}

\textit{Proof}: {One only needs to }apply Theorem \ref{key-thm}.\vspace
{-0.6cm} \begin{flushright}
\mbox{$\Box$}
\end{flushright}\noindent

\subsection{Computations and examples}

We now show how to compute $g(Z)  :=E(\langle DZ,-DL^{-1}%
Z\rangle_{\EuFrak H}|Z)$ in practice. We then provide several examples using
this computation.

\begin{prop}
\label{computation} Write $DZ=\Phi_{Z}(X)$ with a measurable function
$\Phi_{Z}:\mathbb{R}^{\EuFrak H}\rightarrow\EuFrak H$. We have
\[
g(Z)  =\int_{0}^{\infty}e^{-u}\,\mathbf{E}\big(\langle\Phi
_{Z}(X),\Phi_{Z}(e^{-u}X+\sqrt{1-e^{-2u}}X^{\prime})\rangle_{\EuFrak
H}|Z\big)du,
\]
where $X^{\prime}$ stands for an independent copy of $X$, and is such that $X$
and $X^{\prime}$ are defined on the product probability space $(\Omega
\times\Omega^{\prime},\mathscr{F}\otimes\mathscr{F}^{\prime},P\times
P^{\prime})$. Here $\mathbf{E}$ denotes the mathematical expectation with
respect to $P\times P^{\prime}$.
\end{prop}

\textit{Proof}: We follow the arguments contained in Nourdin and Peccati
\cite[Remark 3.6]{NP}. Without loss of generality, we can assume that
$\EuFrak
H=L^{2}(A,\mathscr{A},\mu)$ where $(A,\mathscr{A})$ is a
measurable space and $\mu$ is a $\sigma$-finite measure without atoms. Let us
consider the chaos expansion of $Z$, given by $Z=\sum_{m=1}^{\infty}%
I_{m}(f_{m}),\,\mbox{with $f_m\in\HH^{\odot m}$}.$ Therefore $-L^{-1}%
Z=\sum_{m=1}^{\infty}\frac{1}{m}I_{m}(f_{m})$ and
\[
-D_{a}L^{-1}Z=\sum_{m=1}^{\infty}I_{m-1}(f_{m}(\cdot,a)),\quad a\in A.
\]
On the other hand, we have $D_{a}Z=\sum_{m=1}^{\infty}mI_{m-1}(f_{m}%
(\cdot,a))$. Thus
\begin{align*}
\int_{0}^{\infty}e^{-u}T_{u}(D_{a}Z)du  &  =\int_{0}^{\infty}e^{-u}\left(
\sum_{m=1}^{\infty}me^{-(m-1)u}I_{m-1}(f_{m}(\cdot,a))\right)  du\\
&  =\sum_{m=1}^{\infty}I_{m-1}(f_{m}(\cdot,a)).
\end{align*}
Consequently,
\[
-DL^{-1}Z=\int_{0}^{\infty}e^{-u}T_{u}(DZ)du.
\]
By Mehler's formula (\ref{mehler}), and since $DZ=\Phi_{Z}(X)$ by assumption,
we deduce that
\[
-DL^{-1}Z=\int_{0}^{\infty}e^{-u}E^{\prime}\big(\Phi_{Z}(e^{-u}X+\sqrt
{1-e^{-2u}}X^{\prime})\big)du.
\]
Using $E(E^{\prime}(\ldots)|Z)=\mathbf{E}(\ldots|Z)$, the desired conclusion
follows. \vspace{-0.6cm} \begin{flushright}
\mbox{$\Box$}
\end{flushright}\noindent

By combining (\ref{ohquellebelleformulebis}) with Proposition
\ref{computation}, we get the formula (\ref{rho2}) given in the introduction,
more precisely:

\begin{cor}
\label{main-thm} Let $Z\in\mathbb{D}^{1,2}$ be centered, and let $\Phi
_{Z}:\mathbb{R}^{\EuFrak H}\rightarrow\EuFrak H$ be measurable and such that
$DZ=\Phi_{Z}(X)$. Assume that condition (\ref{sigmamin}) holds. Then $Z$ has a
density $\rho$ given, for almost all $z\in\mathbb{R}$, by
\begin{align*}
\rho(z)  &  =\frac{E|Z|}{2\int_{0}^{\infty}e^{-u}\,\mathbf{E}\big(\langle
\Phi_{Z}(X),\Phi_{Z}(e^{-u}X+\sqrt{1-e^{-2u}}X^{\prime})\rangle_{\EuFrak H}%
|Z=z\big)du}\\
&  \quad\times\mathrm{exp}\left(  -\int_{0}^{z}\frac{x\,dx}{\int_{0}^{\infty
}e^{-u}\,\mathbf{E}\big(\langle\Phi_{Z}(X),\Phi_{Z}(e^{-u}X+\sqrt{1-e^{-2u}%
}X^{\prime})\rangle_{\EuFrak H}|Z=x\big)du}\right)  .
\end{align*}

\end{cor}

\vskip0.5cm

Now, we give several examples of application of this corollary.

\subsubsection{First example: monotone Gaussian functional, finite case.}

Let $N\sim\mathscr{N}_{n}(0,K)$ with $K\in\mathscr{S}_{n}^{+}(\mathbb{R})$,
and $f:\mathbb{R}^{n}\rightarrow\mathbb{R}$ be a $\mathscr{C}^{1}$ function
having bounded derivatives. We assume, without loss of generality, that each
$N_{i}$ has the form $X(h_{i})$, for a certain centered isonormal process $X$
(over some Hilbert space $\EuFrak H$) and certain functions $h_{i}\in\EuFrak
H$. Set $Z=f(N)-E(f(N))$. The chain rule (\ref{chainrule-lipschitz}) implies
that $Z\in\mathbb{D}^{1,2}$ and that $DZ=\Phi_{Z}(N)=\sum_{i=1}^{n}%
\frac{\partial f}{\partial x_{i}}(N)h_{i}$. Therefore
\[
\langle\Phi_{Z}(X),\Phi_{Z}(e^{-u}X+\sqrt{1-e^{-2u}}X^{\prime})\rangle
_{\EuFrak H}=\sum_{i,j=1}^{n}K_{ij}\frac{\partial f}{\partial x_{i}}%
(N)\frac{\partial f}{\partial x_{j}}(e^{-u}N+\sqrt{1-e^{-2u}}N^{\prime}).
\]
(Compare with Lemma 5.3 in Chatterjee \cite{Chatterbox2}). In particular,
Corollary \ref{main-thm} yields the following.

\begin{prop}
\label{fonctionnel}Let $N\sim\mathscr{N}_{n}(0,K)$ with $K\in\mathscr{S}_{n}%
^{+}(\mathbb{R})$, and $f:\mathbb{R}^{n}\rightarrow\mathbb{R}$ be a
$\mathscr{C}^{1}$ function with bounded derivatives. If there exist
$\alpha_{i},\beta_{i}\geqslant0$ such that $\alpha_{i}\leqslant\frac{\partial
f}{\partial x_{i}}(x)\leqslant\beta_{i}$ for any $i\in\{1,\ldots,n\}$ and
$x\in\mathbb{R}^{n}$, if $K_{ij}\geqslant0$ for any $i,j\in\{1,\ldots,n\}$ and
if $\sum_{i,j=1}^{n}\alpha_{i}\alpha_{j}K_{ij}>0$, then $Z=f(N)-E(f(N))$ has a
density $\rho$ satisfying, for almost all $z\in\mathbb{R}$,
\begin{align*}
&  \frac{E|Z|}{2\sum_{i,j=1}^{n}\beta_{i}\beta_{j}K_{ij}}\,\mathrm{exp}\left(
-\frac{z^{2}}{2\sum_{i,j=1}^{n}\alpha_{i}\alpha_{j}K_{ij}}\right)  \\
&  \hskip2cm\leqslant\rho(z)\leqslant\frac{E|Z|}{2\sum_{i,j=1}^{n}\alpha
_{i}\alpha_{j}K_{ij}}\,\mathrm{exp}\left(  -\frac{z^{2}}{2\sum_{i,j=1}%
^{n}\beta_{i}\beta_{j}K_{ij}}\right)  .
\end{align*}

\end{prop}

\vskip.5cm

\subsubsection{Second example: proof of Proposition \ref{rhointroprop}.}

Let $N\sim\mathscr{N}_{n}(0,K)$ with $K\in\mathscr{S}_{n}^{+}(\mathbb{R})$.
Once again, we assume that each $N_{i}$ has the form $X(h_{i})$, for a certain
centered isonormal process $X$ (over some Hilbert space $\EuFrak H$) and
certain functions $h_{i}\in\EuFrak H$. Let $Z=\max N_{i}-E(\max N_{i})$, and
set
\[
I_{u}=\mathrm{argmax}_{1\leqslant i\leqslant n}(e^{-u}X(h_{i})+\sqrt
{1-e^{-2u}}X^{\prime}(h_{i}))\quad\mbox{for $u\geq 0$.}
\]

\begin{lemme}
\label{hi0} For any $u\geqslant0$, $I_{u}$ is a well-defined random element of
$\{1,\ldots,n\}$. Moreover, $Z\in\mathbb{D}^{1,2}$ and we have $DZ=\Phi
_{Z}(N)=h_{I_{0}}$.
\end{lemme}

\textit{Proof}: Fix $u\geqslant0$. Since, for any $i\neq j$, we have
\begin{align*}
&  P\big(e^{-u}X(h_{i})+\sqrt{1-e^{-2u}}X^{\prime}(h_{i})=e^{-u}X(h_{j}%
)+\sqrt{1-e^{-2u}}X^{\prime}(h_{j})\big)\\
&  =P\big(X(h_{i})=X(h_{j})\big)=0,
\end{align*}
the random variable $I_{u}$ is a well-defined element of $\{1,\ldots,n\}$.
Now, if $\Delta_{i}$ denotes the set $\{x\in\mathbb{R}^{n}:x_{j}\leqslant
x_{i}\,\mbox{ for all $j$}\}$, observe that $\frac{\partial}{\partial x_{i}%
}\max=\mathbf{1}_{\Delta_{i}}$ almost everywhere. 
%%J'ai enlevé la phrase ci-dessous
%%Moreover, the function $\max:\mathbb{R}^{n}\rightarrow\mathbb{R}$ is Lipshitz, see (\ref{maxlip}) below.
The desired conclusion follows from the Lipshitz version of the chain rule
(\ref{chainrule-lipschitz}), and the following Lipshitz property of the $\max$
function, which is easily proved by induction on $n\geqslant1$:%
\begin{equation}
\big|\max(y_{1},\ldots,y_{n})-\max(x_{1},\ldots,x_{n})\big|\leqslant\sum
_{i=1}^{n}|y_{i}-x_{i}|\quad\mbox{for any $x,y\in\R^n$.} \label{maxlip}%
\end{equation}
\vspace{-0.6cm} \begin{flushright}
\mbox{$\Box$}
\end{flushright}\noindent In particular, we deduce from Lemma \ref{hi0} that
\begin{equation}
\langle\Phi_{Z}(X),\Phi_{Z}(e^{-u}X+\sqrt{1-e^{-2u}}X^{\prime})\rangle
_{\EuFrak H}=K_{I_{0},I_{u}}. \label{l-1max}%
\end{equation}
By combining this fact with Corollary \ref{main-thm}, we get Proposition
\ref{rhointroprop}, which we restate.

\begin{prop}
\label{finite-max} Let $N\sim\mathscr{N}_{n}(0,K)$ with $K\in\mathscr{S}_{n}%
^{+}(\mathbb{R})$. If there exists $\sigma_{\min},\sigma_{\max}>0$ such that
$\sigma_{\min}^{2}\leqslant K_{ij}\leqslant\sigma_{\max}^{2}$ for any
$i,j\in\{1,\ldots,n\}$, then $Z=\max N_{i}-E(\max N_{i})$ has a density $\rho$
satisfying (\ref{rho-intro}) for almost all $z\in\mathbb{R}$.
\end{prop}

\vskip.5cm

\subsubsection{Third example: monotone Gaussian functional, continuous case.}

Assume that $X=(X_{t},\,t\in\lbrack0,T])$ is a centered Gaussian process with
continuous paths, and that $f:\mathbb{R}\rightarrow\mathbb{R}$ is
$\mathscr{C}^{1}$ with a bounded derivative. Consider $Z=\int_{0}^{T}%
f(X_{v})dv-E\left(  \int_{0}^{T}f(X_{v})dv\right)  $. Then $Z\in
\mathbb{D}^{1,2}$ and we have $DZ=\Phi_{Z}(X)=\int_{0}^{T}f^{\prime}%
(X_{v})\mathbf{1}_{[0,v]}dv$. Therefore
\begin{align*}
&  \langle\Phi_{Z}(X),\Phi_{Z}(e^{-u}X+\sqrt{1-e^{-2u}}X^{\prime}%
)\rangle_{\EuFrak H}\\
&  =\iint_{[0,T]^{2}}f^{\prime}(X_{v})f^{\prime}(e^{-u}X_{w}+\sqrt{1-e^{-2u}%
}X_{w}^{\prime})E(X_{v}X_{w})dvdw.
\end{align*}
Using Corollary \ref{main-thm}, we get the following.

\begin{prop}
\label{infdimfonc}Assume that $X=(X_{t},\,t\in\lbrack0,T])$ is a centered
Gaussian process with continuous paths, and that $f:\mathbb{R}\rightarrow
\mathbb{R}$ is $\mathscr{C}^{1}$. If there exists $\alpha,\beta,\sigma_{\min
},\sigma_{\max}>0$ such that $\alpha\leqslant f^{\prime}(x)\leqslant\beta$ for
all $x\in\mathbb{R}$ and $\sigma_{\min}^{2}\leqslant E(X_{v}X_{w}%
)\leqslant\sigma_{\max}^{2}$ for all $v,w\in\lbrack0,T]$, then $Z=\int_{0}%
^{T}f(X_{v})dv-E\left(  \int_{0}^{T}f(X_{v})dv\right)  $ has a density $\rho$
satisfying, for almost all $z\in\mathbb{R}$,
\[
\frac{E|Z|}{2\beta^{2}\,\sigma_{\max}^{2}\,T^{2}}\,e^{-\frac{z^{2}}%
{2\alpha^{2}\,\sigma_{\min}^{2}T^{2}}}\leqslant\rho(z)\leqslant\frac
{E|Z|}{2\alpha^{2}\,\sigma_{\min}^{2}\,T^{2}}\,e^{-\frac{z^{2}}{2\beta
^{2}\,\sigma_{\max}^{2}T^{2}}}.
\]

\end{prop}

\vskip.5cm

\subsubsection{Fourth example: supremum of a Gaussian process}

Fix $a<b$, and assume that $X=(X_{t},\,t\in\lbrack a,b])$ is a centered
Gaussian process with continuous paths and such that $E|X_{t}-X_{s}|^{2}\neq0$
for all $s\neq t$. Set $Z=\sup_{[a,b]}X-E(\sup_{[a,b]}X)$, and let $\tau_{u}$
be the (unique) random point where $e^{-u}X+\sqrt{1-e^{-2u}}X^{\prime}$
attains its maximum on $[a,b]$. Note that $\tau_{u}$ is well-defined, see e.g.
Lemma 2.6 in \cite{kimpollard}. Moreover, we have that $Z\in\mathbb{D}^{1,2}$,
see Proposition 2.1.10 in \cite{nualartbook}, and $DZ=\Phi_{Z}(X)=\mathbf{1}%
_{[0,\tau_{0}]}$, see Lemma 3.1 in \cite{decreununu}. Therefore
\[
\langle\Phi_{Z}(X),\Phi_{Z}(e^{-u}X+\sqrt{1-e^{-2u}}X^{\prime})\rangle
_{\EuFrak H}=R(\tau_{0},\tau_{u})
\]
where $R(s,t)=E(X_{s}X_{t})$ is the covariance function of $X$. Using
Corollary \ref{main-thm}, the following obtains.

\begin{prop}
\label{max-intro} Let $X=(X_{t},\,t\in\lbrack a,b])$ be a centered Gaussian
process with continuous paths, and such that $E|X_{t}-X_{s}|^{2}\neq0$ for all
$s\neq t$. Assume that, for some real $\sigma_{\min},\sigma_{\max}>0$, we have
$\sigma_{\min}^{2}\leqslant E(X_{s}X_{t})\leqslant\sigma_{\max}^{2}$ for any
$s,t\in\lbrack a,b]$. Then, $Z=\sup_{[a,b]}X-E(\sup_{[a,b]}X)$ has a density
$\rho$ satisfying, for almost all $z\in\mathbb{R}$,
\[
\frac{E|Z|}{2\sigma_{\max}^{2}}\,e^{-\frac{z^{2}}{2\,\sigma_{\min}^{2}}%
}\leqslant\rho(z)\leqslant\frac{E|Z|}{2\sigma_{\min}^{2}}\,e^{-\frac{z^{2}%
}{2\,\sigma_{\max}^{2}}}.
\]

\end{prop}

To the best of our knowledge, Proposition \ref{max-intro}, as well as
Proposition \ref{finite-max}, contain the first bounds ever established for
the \textit{density} of the supremum of a general Gaussian process. When
integrated over $z$, the upper bound above improves the classical
concentration inequalities (\ref{vitale-b}), (\ref{borell-b}),
(\ref{vitale-borell-b}) on the tail of $Z$, see e.g. the upper bound in
(\ref{improve-intro}); the lower bound for the left-hand tail of $Z$ which one
obtains by integration, appears to be entirely new. When applied to the case
of fractional Brownian motion, we get the following.

\begin{cor}
\label{max-fbm} Let $b>a>0$, and $B=(B_{t},\,t\geqslant0)$ be a fractional
Brownian motion with Hurst index $H\in\lbrack1/2,1)$. Then the random variable
$Z=\sup_{[a,b]}B-E\big(\sup_{[a,b]}B\big)$ has a density $\rho$ satisfying
(\ref{ffbbmm}) for almost all $z\in\mathbb{R}$.
\end{cor}

\textit{Proof}: For any choice of the Hurst parameter $H\in(1/2,1)$, the
Gaussian space generated by $B$ can be identified with an isonormal Gaussian
process of the type $X=\{X(h):h\in\EuFrak H\}$, where the real and separable
Hilbert space $\EuFrak H$ is defined as follows: (i) denote by $\mathscr{E}$
the set of all $\mathbb{R}$-valued step functions on $\mathbb{R}_{+}$, (ii)
define $\EuFrak H$ as the Hilbert space obtained by closing $\mathscr{E}$ with
respect to the scalar product
\[
\left\langle {\mathbf{1}}_{[0,t]},{\mathbf{1}}_{[0,s]}\right\rangle
_{\EuFrak H}=E(B_{t}B_{s})=\frac12\big(t^{2H}+s^{2H}-|t-s|^{2H}\big).
\]
In particular, with such a notation, one has that $B_{t}=X(\mathbf{1}%
_{[0,t]})$. The reader is referred e.g. to \cite{nualartbook} for more details
on fractional Brownian motion.

Now, the desired conclusion is a direct application of Proposition
\ref{max-intro} since, for all $a\leqslant s<t\leqslant b$,
\[
E(B_{s}B_{t})\leqslant\sqrt{E(B_{s}^{2})}\,\sqrt{E(B_{t}^{2})}=(st)^{H}%
\leqslant b^{2H}%
\]
and
\begin{align*}
E(B_{s}B_{t})  &  =\frac{1}{2}\big(t^{2H}+s^{2H}-(t-s)^{2H}\big)=H(2H-1)\iint
_{[0,s]\times\lbrack0,t]}|v-u|^{2H-2}dudv\\
&  \geqslant H(2H-1)\iint_{[0,a]\times\lbrack0,a]}|v-u|^{2H-2}dudv=E(B_{a}%
^{2})=a^{2H}.
\end{align*}
\vspace{-0.6cm} \begin{flushright}
\mbox{$\Box$}
\end{flushright}\noindent\vspace{0.1in}

\section{Concentration inequalities}

\label{sect-concentration} Now, we investigate what can be said when $g(Z)  
=E(\langle DZ,-DL^{-1}Z\rangle_{\EuFrak H}|Z)$ just admits a lower
(resp. upper) bound. Results under such hypotheses are more difficult to
obtain than in the previous section, since there we could use bounds on
$g(Z)  $ in both directions to good effect; this is apparent, for
instance, in the appearance of both the lower and upper bounding values
$\sigma_{\min}$ and $\sigma_{\max}$ in each of the two bound in
(\ref{rho-intro}), or more generally in Corollary \ref{general}. However,
given our previous work, tails bounds can be readily obtained: most of the
analysis of the role of $g(Z)  $ in tail estimates is already
contained in the proof of Theorem \ref{key-thm}.

\subsection{Upper bounds}

{Our first result allows comparisons both to the Gaussian and exponential
tails}.

\begin{thm}
\label{concentration} Let $Z\in\mathbb{D}^{1,2}$ with zero mean, 
$g(Z)  =E(\langle DZ,-DL^{-1}Z\rangle_{\EuFrak H}|Z)$, and fix
$\alpha\geqslant0$ and $\beta>0$. Assume that

\begin{enumerate}
\item[(i)] $g(Z)  \leqslant\alpha Z+\beta$ almost surely;

\item[(ii)] $Z$ has a density $\rho$.
\end{enumerate}

Then, for all $z>0$, we have
\[
P(Z\geqslant z)\leqslant\mathrm{exp}\left(  -\frac{z^{2}}{2\alpha z+2\beta
}\right)  .
\]

\end{thm}

\textit{Proof}: We follow the same line of reasoning as in \cite[Theorem
1.5]{Chatterbox}. For any $A>0$, define $m_{A}:[0,+\infty)\rightarrow
\mathbb{R}$ by $m_{A}(\theta)=E\big(e^{\theta Z}\mathbf{1}_{\{Z\leqslant
A\}}\big)$. By Lebesgue differentiation theorem, we have
\[
m_{A}^{\prime}(\theta)=E(Ze^{\theta Z}\mathbf{1}_{\{Z\leqslant A\}}%
)\quad\mbox{for all $\theta\geq 0$}.
\]
Therefore, we can write
\begin{align*}
m_{A}^{\prime}(\theta)  &  =\int_{-\infty}^{A}z\,e^{\theta z}\,\rho(z)dz\\
&  =-e^{\theta A}\int_{A}^{\infty}y\rho(y)dy+\theta\int_{-\infty}^{A}e^{\theta
z}\left(  \int_{z}^{\infty}y\rho(y)dy\right)  dz\quad
\mbox{by integration by parts}\\
&  \leqslant\theta\int_{-\infty}^{A}e^{\theta z}\left(  \int_{z}^{\infty}%
y\rho(y)dy\right)  dz\quad\mbox{{since $\int_A^\infty
y\rho(y)dy\geq 0$}}\\
&  =\theta E\big(g(Z)  \,e^{\theta Z}\,\mathbf{1}_{\{Z\leqslant
A\}}\big),
\end{align*}
where the last line follows from identity (\ref{ohlabelleformule}). Due to the
assumption $(i)$, we get
\[
m_{A}^{\prime}(\theta)\leqslant\theta\,\alpha\,m_{A}^{\prime}(\theta
)+\theta\,\beta\,m_{A}(\theta),
\]
that is, for any $\theta\in(0,1/\alpha)$:
\[
\frac{m_{A}^{\prime}(\theta)}{m_{A}(\theta)}\leqslant\frac{\theta\beta
}{1-\theta\alpha}.
\]
By integration and since $m_{A}(0)=P(Z\leqslant A)\leqslant1$, this gives, for
any $\theta\in(0,1/\alpha)$:
\[
m_{A}(\theta)\leqslant\mathrm{exp}\left(  \int_{0}^{\theta}\frac{\beta
u}{1-\alpha u}du\right)  \leqslant\mathrm{exp}\left(  \frac{\beta\theta^{2}%
}{2(1-\theta\alpha)}\right)  .
\]
Using Fatou's inequality (as $A\rightarrow\infty$) in the previous relation
implies
\[
E\big(e^{\theta Z}\big)\leqslant\mathrm{exp}\left(  \frac{\beta\theta^{2}%
}{2(1-\theta\alpha)}\right)
\]
for all $\theta\in(0,1/\alpha)$. Therefore, for all $\theta\in(0,1/\alpha)$,
we have
\[
P(Z\geqslant z)=P(e^{\theta Z}\geqslant e^{\theta z})\leqslant e^{-\theta
z}E\big(e^{\theta Z}\big)\leqslant\mathrm{exp}\left(  \frac{\beta\theta^{2}%
}{2(1-\theta\alpha)}-\theta z\right)  .
\]
Choosing $\theta=\frac{z}{\alpha z+\beta}\in(0,1/\alpha)$ gives the desired
result. \vspace{-0.6cm} \begin{flushright}
\mbox{$\Box$}
\end{flushright}\noindent

Let us give an example of application of Theorem \ref{concentration}. {{Assume
that $B=(B_{t},\,t\geqslant0)$ is a fractional Brownian motion with Hurst
index $H\in(0,1)$. For any choice of the parameter $H$, as already mentioned
in the proof of Corollary \ref{max-fbm}, the Gaussian space generated by $B$
can be identified with an isonormal Gaussian process of the type
$X=\{X(h):h\in\EuFrak H\}$, where the real and separable Hilbert space
$\EuFrak H$ is defined as follows: (i) denote by $\mathscr{E}$ the set of all
$\mathbb{R}$-valued step functions on $\mathbb{R}_{+}$, (ii) define $\EuFrak
H$ as the Hilbert space obtained by closing $\mathscr{E}$ with respect to the
scalar product
\[
\left\langle {\mathbf{1}}_{[0,t]},{\mathbf{1}}_{[0,s]}\right\rangle
_{\EuFrak
H}=E(B_{t}B_{s})=\frac{1}{2}\big(t^{2H}+s^{2H}-|t-s|^{2H}\big).
\]
In particular, with such a notation one has that $B_{t}=X(\mathbf{1}_{[0,t]}%
)$. Now, }}let
\[
Z=Z_{T}:=\int_{0}^{T}B_{u}^{2}du-\frac{T^{2H+1}}{2H+1}.
\]
By the scaling property of fractional Brownian motion, we see first that
$Z_{T}$ has the same distribution as $T^{2H+1}Z_{1}$. Thus we choose $T=1$
without loss of generality; we denote $Z=Z_{1}$. Now observe that
$Z\in\mathbb{D}^{1,2}$ lives in the second Wiener chaos of $B$. In particular,
we have $-L^{-1}Z=\frac{1}{2}Z$. Moreover $DZ=2\int_{0}^{1}B_{u}%
\,\mathbf{1}_{[0,u]}du,$ so that
\begin{align*}
\langle DZ,-DL^{-1}Z\rangle_{\EuFrak H}  &  =\frac{1}{2}\Vert DZ\Vert
_{\EuFrak H}^{2}=2\int_{[0,1]^{2}}B_{u}B_{v}\,E(B_{u}B_{v})dudv\\
&  \leqslant2\int_{[0,1]^{2}}|B_{u}|\,|B_{v}|\,\big|E(B_{u}B_{v})\big|dudv\\
&  \leqslant2\int_{[0,1]^{2}}|B_{u}|\,|B_{v}|\,u^{H}v^{H}dudv=2\left(
\int_{0}^{1}|B_{u}|u^{H}du\right)  ^{2}\\
&  \leqslant2\int_{0}^{1}B_{u}^{2}du\times\int_{0}^{1}u^{2H}du=\frac{1}%
{H+1/2}\int_{0}^{1}B_{u}^{2}du\\
&  =\frac{1}{H+1/2}\left(  Z+\frac{1}{2H+1}\right)  .
\end{align*}
Since it is easily shown that $Z$ has a density, Theorem \ref{concentration}
implies the desired conclusion in Proposition \ref{fbm}, or with $c_{H}%
=H+1/2$,%
\[
P\left(  Z_{1}\geqslant z\right)  \leqslant\exp\left(  -\frac{z^{2}c_{H}^{2}%
}{2c_{H}z+1}\right)  .
\]
{By scaling}, this shows that the tail of $Z_{T}/T^{2H+1}$ behaves
asymptotically like that of an exponential random variable with mean
$\nu=\left(  H/2+1/4\right)  ^{-1}$.

For the moment, it is not possible to use our tools to investigate a lower
bound on this tail, see the forthcoming Section \ref{lowerb}. We have also
investigated the possibility of using such tools as the formula (\ref{rho}),
or the density lower bounds found in \cite{Ncras}, thinking that a specific
second-chaos situation might be tractable despite the reliance on the
divergence operator, but these tools seem even less appropriate. However, in
this particular instance, we can perform a calculation by hand, as follows. By
Jensen's inequality, with $\mu=\left(  2H+1\right)  ^{-1}$, we have that
$Z+\mu=Z_{1}+\mu\geqslant\left(  \int_{0}^{1}B_{u}du\right)  ^{2}$. Thus%
\[
P\big(Z_{1}\geqslant z\big)\geqslant P\left(  \left(  \int_{0}^{1}%
B_{u}du\right)  ^{2}\geqslant z+\mu\right)  =P\left(  \left\vert \int_{0}%
^{1}B_{u}du\right\vert \geqslant\sqrt{z+\mu}\right)  .
\]
Here of course, the random variable $N=\int_{0}^{1}B_{u}du$ is centered
Gaussian, and its variance can be calculated by hand:
\begin{align*}
\sigma^{2}:=E\left(  N^{2}\right)   &  =\iint_{[0,1]^{2}}E\left(  B_{u}%
B_{v}\right)  dudv\\
&  =\int_{0}^{1}dv\int_{0}^{v}du\left(  u^{2H}+v^{2H}-\left(  v-u\right)
^{2H}\right)  =\frac{1}{2H+2}.
\end{align*}
Therefore, by the standard lower bound on the tail of a Gaussian r.v., {{that
is $\int_{z}^{\infty}e^{-y^{2}/2}dy\geqslant\frac{z}{1+z^{2}}e^{-z^{2}/2}$ for
all $z>0$, we get}}
\begin{align}
P\left(  Z_{1}\geqslant z\right)   &  \geqslant\frac{{{\sigma\sqrt{z+\mu}}}%
}{\sigma^{2}+z+\mu}\exp\left(  -\frac{z+\mu}{2\sigma^{2}}\right) \nonumber\\
&  {{\underset{z\rightarrow\infty}{\sim}}}\frac{1}{{{\sqrt{z}\sqrt{2H+2}}}%
}\exp\left(  -\frac{H+1}{2H+1}\right)  \exp\left(  -\left(  H+1\right)
z\right)  . \label{lbfBm}%
\end{align}
%where the last relation is asymptotic equivalence for large $z$.
Abusively ignoring the factor ${{z^{-1/2}}}$ in this lower bound, we can
summarize our results by saying that $Z_{T}/T^{2H+1}$ has a tail that is
bounded above and below by exponential tails with respective means $\left(
H/2+1/4\right)  ^{-1}$ and $\left(  H+1\right)  ^{-1}$.\bigskip

As another example, let us explain how Theorem \ref{concentration} allows to
easily recover both the Borell-Sudakov-type inequalities (\ref{vitale-b}) and
(\ref{borell-b}), for $Z$ defined as the {centered} supremum of a Gaussian
vector in (\ref{def-z}). We can assume, without loss of generality, that each
$N_{i}$ has the form $X(h_{i})$, for a certain centered isonormal process $X$
(over some Hilbert space $\EuFrak H$) and certain functions $h_{i}\in\EuFrak
H$. Condition $(ii)$ of Theorem \ref{concentration} is easily satisfied while
for condition $(i)$, we have, by combining (\ref{l-1max}) with Proposition
\ref{computation}:
\begin{equation}
\langle DZ,-DL^{-1}Z\rangle_{\EuFrak H}=\int_{0}^{\infty}e^{-u}K_{I_{0},I_{u}%
}du\leqslant\max_{1\leqslant i,j\leqslant n}K_{ij}\,=\,\sigma_{\max}^{2}
\label{rty}%
\end{equation}
so that
\[
g(Z)  \leqslant\sigma_{\max}^{2}\quad\mbox{almost surely}.
\]
In other words, condition $(i)$ is satisfied with $\alpha=0$ and $\beta
=\sigma_{\max}^{2}$. Therefore $P(Z\geqslant z)\leqslant\mathrm{exp}\left(
-\frac{z^{2}}{2\,\sigma_{\max}^{2}}\right)  $, for all $z>0$, and
(\ref{vitale-b}) is shown. The proof of (\ref{borell-b}) follows the same
lines, by considering $-Z$ instead of $Z$.

\subsection{Lower bounds}

\label{lowerb}

We now investigate a lower bound analogue of Theorem \ref{concentration}.
Recall we still use the notation $g(z) =E(\langle
DZ,-DL^{-1}Z\rangle_{\EuFrak H}|Z=z)$.

\begin{thm}
\label{thmlbnostein} Let $Z\in\mathbb{D}^{1,2}$ with zero mean, and fix
$\sigma_{\min},\alpha>0$ and $\beta>1$. Assume that

\begin{enumerate}
\item[(i)] $g(Z)  \geqslant\sigma_{\min}^{2}$ almost surely.
\end{enumerate}

The existence of the density $\rho$ of $Z$ is thus ensured by Theorem
\ref{key-thm}. Also assume that

\begin{enumerate}
\item[(ii)] the function $h\left(  x\right)  :=x^{1+\beta}\rho\left(
x\right)  \ $is decreasing on $[\alpha,+\infty)$.
\end{enumerate}

Then, for all $z\geqslant\alpha$, we have
\[
P(Z\geqslant z)\geqslant\frac{1}{2}\left(  1-\frac{1}{\beta}\right)
E|Z|~\frac{1}{z}\,\mathrm{exp}\left(  -\frac{z^{2}}{2\,\sigma_{\min}^{2}%
}\right)  .
\]

Alternately, {instead of $(ii)$}, assume that there exists $0<\alpha<2$ such that

\begin{enumerate}
\item[(ii)'] $\limsup_{z\rightarrow\infty}z^{-\alpha}\log g(z)
<\infty$.
\end{enumerate}

Then, for any $\varepsilon>0$, there exist $K,z_{0}>0$ such that, for all
$z>z_{0}$,%
\[
P(Z\geqslant z)\geqslant K~\mathrm{exp}\left(  -\frac{z^{2}}{\left(
2-\varepsilon\right)  \sigma_{\min}^{2}}\right)  .
\]

\end{thm}

\textit{Proof}: {{First, let us relate the function $\varphi(z)=\int
_{z}^{\infty}y\rho(y)dy$ to the tail of $Z$. By integration by parts, we get
\begin{equation}
\varphi\left(  z\right)  =z\,P(Z\geqslant z)+\int_{z}^{\infty}P(Z\geqslant
y)dy.\label{r}%
\end{equation}
}}If we assume $(ii)$, since $h$ is decreasing, for any $y>z\geqslant\alpha$
we have $\frac{y\rho\left(  y\right)  }{z\rho(z)}\leqslant\left(  \frac{z}%
{y}\right)  ^{\beta}.$ Then we have, for any $z\geqslant\alpha$:
\[
P(Z\geqslant z)=z\rho\left(  z\right)  \int_{z}^{\infty}\frac{1}{y}\frac
{y\rho\left(  y\right)  }{z\rho\left(  z\right)  }dy\leqslant z\rho\left(
z\right)  z^{\beta}\int_{z}^{\infty}\frac{dy}{y^{1+\beta}}=\frac{z\rho\left(
z\right)  }{\beta}.
\]
By putting that inequality into (\ref{r}), we get%
\[
\varphi(z)\leqslant z\,P(Z\geqslant z)+\frac{1}{\beta}\int_{z}^{\infty}%
y\rho(y)dy=z\,P(Z\geqslant z)+\frac{1}{\beta}\,\varphi(z)
\]
so that $P(Z\geqslant z)\geqslant\left(  1-\frac{1}{\beta}\right)
\frac{\varphi(z)}{z}.$ Combined with (\ref{Flb}), this gives the desired conclusion.

Now assume $(ii)^{\prime}$ instead. Here the proof needs to be modified. From
the key result of Theorem \ref{key-thm} and condition $(i)$, we have%
\[
\rho(z)\geqslant\frac{E|Z|}{2\,g(z)  }\,\mathrm{exp}\left(
-\frac{z^{2}}{2\,\sigma_{\min}^{2}}\right)  .
\]
Let $\Psi\left(  z\right)  $ denote the unnormalized Gaussian tail $\int
_{z}^{\infty}\mathrm{exp}\left(  -\frac{y^{2}}{2\,\sigma_{\min}^{2}}\right)
dy$. We can write, using the Schwarz inequality,
\begin{align*}
\Psi^{2}(z) &  =\left(  \int_{z}^{\infty}\mathrm{exp}\left(  -\frac{y^{2}%
}{2\,\sigma_{\min}^{2}}\right)  \sqrt{g(y)}~\frac{1}{\sqrt{g(y)}}~dy\right)
^{2}\\
&  \leqslant\int_{z}^{\infty}\mathrm{exp}\left(  -\frac{y^{2}}{2\,\sigma
_{\min}^{2}}\right)  ~g(y)~dy\times\int_{z}^{\infty}\mathrm{exp}\left(
-\frac{y^{2}}{2\,\sigma_{\min}^{2}}\right)  ~\frac{1}{g(y)}~dy
\end{align*}
so that
\begin{align*}
P(Z\geqslant z) &  =\int_{z}^{\infty}\rho\left(  y\right)  dy\\
&  \geqslant\frac{E\left\vert Z\right\vert }{2}\int_{z}^{\infty}%
e^{-y^{2}/\left(  2\sigma_{\min}^{2}\right)  }~\frac{1}{g(y)
}~dy\\
&  \geqslant\frac{E\left\vert Z\right\vert }{2}\frac{\Psi^{2}\left(  z\right)
}{\int_{z}^{\infty}e^{-y^{2}/\left(  2\sigma_{\min}^{2}\right)  }g\left(
y\right)  dy}.
\end{align*}
Using the classical inequality $\int_{z}^{\infty}e^{-y^{2}/2}dy\geqslant
\frac{z}{1+z^{2}}e^{-z^{2}/2}$, we get
\begin{equation}
P(Z\geqslant z)\geqslant\frac{E|Z|}{2}~\frac{\sigma_{\min}^{4}z^{2}%
}{\big(\sigma_{\min}^{2}+z^{2}\big)^{2}}~\frac{\mathrm{exp}\left(
-\frac{z^{2}}{\sigma_{\min}^{2}}\right)  }{\int_{z}^{\infty}\mathrm{exp}%
\left(  -\frac{y^{2}}{2\sigma_{\min}^{2}}\right)  g(y)dy}.\label{iv}%
\end{equation}

Under condition $(ii)^{\prime}$, we have that there exists $c>0$ such that,
for $y$ large enough, $g(y) \leqslant e^{cy^{\alpha}}$ with $0<\alpha<2$. We
leave it to the reader to check that the conclusion now follows by an
elementary calculation from (\ref{iv}). \vspace{-0.6cm} \begin{flushright}
\mbox{$\Box$}
\end{flushright}
\noindent

\begin{rem}

%%J'ai rajouté une numérotation
\begin{enumerate}
\item
Inequality (\ref{iv}) itself may be of independent
interest, when the growth of $g$ can be controlled, but not
as efficiently as in $(ii)^{\prime}$.
\item
Condition $(ii)$ implies that $Z$ has a moment of
order greater than $\beta$. Therefore it can be considered as a
technical regularity and integrability condition. Condition $(ii)^{\prime}$
may be easier to satisfy in cases where a good handle on 
$g$ exists. Yet the use of the Schwarz inequality in the above proof
means that conditions $(ii)^{\prime}$ is presumably stronger than it
needs to be.
\item
In general, one can see that deriving lower bounds on tails of random
variables with little upper bound control is a difficult task, deserving of
further study.
\end{enumerate}
\end{rem}

\vskip1cm

\noindent\textbf{Acknowledgment}: We are grateful to Paul Malliavin, David
Nualart, and Giovanni Peccati, for helpful comments.

\end{document}